\documentclass[a4paper,oneside,10pt]{article}
\usepackage{amssymb}
\usepackage{amsmath}
\usepackage{amsthm}
\usepackage{fullpage}
\usepackage[cp1250]{inputenc}
\usepackage[IL2]{fontenc}
\newtheorem{Theorem}{Theorem}
\newtheorem{Lemma}[Theorem]{Lemma}

\begin{document}
\author{Jan Hladk\'y
\thanks{Department of Applied Mathematics, Faculty of Mathematics and Physics, Charles University, Malostransk\'e n\'am\v est\'i
 25, 118 00, Prague, Czech Republic and Zentrum
 Mathematik, Technische Universit\"at M\"unchen,
 Boltzmannstra\ss{}e~3, D-85747 Garching bei M\"unchen,
 Germany. E-mail: {\tt
 honzahladky@googlemail.com}. Research was supported by
 the Grant Agency of Charles University, and by the DFG Research Training Group 1408 ``Methods for Discrete Structures'' via its 2007 predoc course ``Integer Points in Polyhedra''. }
       \and
        Mathias Schacht
\thanks{Institut f\"ur Informatik, Humboldt-Universit\"at zu Berlin, Unter den Linden 6, D-10099 Berlin, Germany. E-mail: {\tt schacht@informatik.hu-berlin.de}}}
\title{Note on bipartite graph tilings}
\maketitle \begin{abstract}Let $s < t$ be two fixed positive integers. We study sufficient minimum degree conditions for a bipartite
graph $G$, with both color classes of size $n=k(s+t)$, which ensure that
$G$ has a $K_{s,t}$-factor. 
Our result extends the work of Zhao, who determined the
minimum degree threshold which guarantees that a bipartite
graph has a $K_{s,s}$-factor.
\end{abstract}
 
\section{Introduction}
For two (finite, loopless, simple) graphs $H$ and $G$, we
say that $G$ contains an {\em $H$-factor} if there exist
$v(G)/v(H)$ vertex-disjoint copies of $H$ in $G$. As a
synonym, we say that there exists an {\em $H$-tiling of
$G$}. Obviously, if $G$ contains an $H$-factor, then $v(G)$
is a multiple of $v(H)$. For a fixed graph $H$, necessary
and sufficient conditions on the minimum-degree of $G$
which guarantee that $G$ contains an $H$-factor were
studied extensively. Results in this spirit include the
Tutte 1-factor Theorem (see~\cite{LP86}), the
Hajnal-Szemer\'edi Theorem~\cite{hajnal70:_proof_p}, and
series of improving results by Alon and
Yuster~\cite{alon92:_yuster_GC,alon96:_yuster_JCTB},
Koml\'os~\cite{komlosTT}, Zhao and
Shokoufandeh~\cite{Zhao:komlosConj}, and by K\"uhn and
Osthus~\cite{Kuhn:Tiling}. In~\cite{Kuhn:Tiling} the answer
to the problem is settled (up to a constant) for any $H$. 
It was shown that the relevant parameters are the
chromatic number and the critical chromatic number of $H$.

The additional information that $G$ is $r$-partite might
help to decrease the minimum-degree threshold for $G$
containing an $H$-factor. 
The conjectured
$r$-partite version of the Hajnal-Szemer\'edi
Theorem~\cite{FischerHS} is  such an example. 
(Recently a proof of the approximate version
of the $r$-partite Hajnal-Szemer\'edi Theorem was announced
by Csaba.) In this paper we determine the threshold for
the minimum-degree of a balanced bipartite graph $G$ which
guarantees that $G$ contains a $K_{s,t}$-factor, for
arbitrary integers $s<t$. If the cardinalities of both
color classes of $G$ are $n$, a necessary
condition for $G$ having a $K_{s,t}$-factor is that $n$ is
a multiple of $s+t$. The sufficient minimum-degree
condition is given in Theorem~\ref{thm_upper}, and a
matching lower bound is provided in
Theorem~\ref{thm_lower}. Our work can be seen as an
extension of the work of Zhao~\cite{ZhaoBipTil}, who
investigated the case $s=t$.
\begin{Theorem}[Zhao,~\cite{ZhaoBipTil}]\label{thm_upperZ}
For each $s\ge 2$ there exists a number $k_0$ such that if $G=(A,B; E)$ is a bipartite graph, 
$|A|=|B|=n=ks$, where $k>k_0$, and
$$\delta(G)\geq\Big\{\begin{array}{ll}
                   \frac{n}{2}+s-1   &\mbox{if $k$ is even,} \\
                   \frac{n+3s}{2}-2 &\mbox{if $k$ is odd,} \\
                 \end{array}
$$
then $G$ has a $K_{s,s}$-factor.
\end{Theorem}
Moreover, Zhao showed that the bounds in Theorem~\ref{thm_upperZ} are tight.
We extend those results to $K_{s,t}$-factors with $s<t$.
\begin{Theorem}\label{thm_upper}
Let $1\leq s<t$ be fixed integers. There exists a number $k_0\in
\mathbb{N}$ such that if $G=(A,B; E)$ is a bipartite graph,
$|A|=|B|=n=k(s+t)$, with $k>k_0$, and
$$\delta(G)\geq\Big\{\begin{array}{ll}
                   \frac{n}{2}+s-1   &\mbox{if $k$ is even,} \\
                   \frac{n+t+s}{2}-1 &\mbox{if $k$ is odd,} \\
                 \end{array}
$$
then $G$ has a $K_{s,t}$-factor.
\end{Theorem}
On the other hand, we show that these bounds are best possible.
\begin{Theorem}\label{thm_lower}
Let $1\leq s<t$ be fixed integers. There exists a number $k_0\in
\mathbb{N}$ such that for every $k>k_0$ there exists a bipartite
graph $G=(A,B; E)$, $|A|=|B|=k(s+t)=n$, such that
$$\delta(G)=\Big\{\begin{array}{ll}
                   \frac{n}{2}+s-2   &\mbox{if $k$ is even,} \\
                   \frac{n+t+s}{2}-2 &\mbox{if $k$ is odd and $t\leq 2s+1$,}\\
                 \end{array}
$$
and $G$ does not have a $K_{s,t}$-factor.
\end{Theorem}
The bounds in Theorem~\ref{thm_upper} and~\ref{thm_lower} exhibit a
somewhat surprising phenomenon: for the case when $k$ is even
the bound is independent of the value $t$, while for the case $k$ is odd,
the minimum-degree condition depends on $t$. Moreover, we note that our results are not tight 
for the case~$t> 2s+1$ and~$k$ odd. We are very grateful to Andrzej Czygrinow and Louis DeBiasio for 
drawing our attention to an oversight in Theorem~\ref{thm_lower} in an 
earlier version of this note.

\section{Lower bound}
In this section we prove Theorem~\ref{thm_lower}. We treat three cases (based on the parity of $k$ 
and on the relation between $s$ and $t$) separately. The proof of Theorem~\ref{thm_lower} is constructive, i.e., we will construct
a graph $G$ with the demanded minimum-degree and then argue that $G$ does not contain a $K_{s,t}$-factor.

The building blocks of our constructions are the graphs $P(m,p)$, where $m,p\in \mathbb{N}$. 
The graphs $P(m,p)$ were introduced in~\cite{ZhaoBipTil}. We just state their properties, which will be used throughout this section. 
\begin{Lemma}\label{lemmaP}
For any $p\in \mathbb{N}$ there exists a number $m_0$ such that for any $m\in\mathbb{N}, m>m_0$ there exists a bipartite graph $P(m,p)=(P_1,P_2;E_P)$  satisfying
\begin{itemize}
 \item $|P_1|=|P_2|=m$,
\item $P(m,p)$ is $p$-regular, and
\item $P(m,p)$ does not contain a copy of $K_{2,2}$.
\end{itemize}
\end{Lemma}

\smallskip In all constructions we assume that $n$ is large enough. 
\subsection{Case $k$ is even}
For two integers $m$ and $q$ we write $Q(m,q)$ to denote (any of
possibly many) bipartite graph $Q(m,q)=(Q_1,Q_2; E_Q)$ with the
following properties:
\begin{itemize}
\item $|Q_1|=m, |Q_2|=m-2$,
\item $Q(m,q)$ does not contain any $K_{2,2}$,
\item $\deg(x)\in\{q-1,q\}$ for any vertex $x\in Q_1$, and
\item $\deg(y)=q$ for any vertex $y\in Q_2$.
\end{itemize}
Such graphs $Q(m,q)$ do exist for fixed $q$ and large $m$. One way to
construct them is by taking the graph $P(m,q)=(P_1, P_2; E_P)$ from Lemma~\ref{lemmaP}, selecting two vertices $w_1, w_2\in P_2$ such that they
do not share a common neighbor in $P_1$, and then take $Q(m,q)$ to be the
subgraph of $P(m,q)$ induced by the vertex sets $P_1$,
$P_2\setminus\{w_1,w_2\}$. In particular, the graph $Q(m,0)$ is the
empty graph.

Now we describe the construction of the graph $G$. Partition
$A=A_1+A_2$, $B=B_1+B_2$, $|A_1|=|B_1|=\frac{n}{2}+1$,
$|A_2|=|B_2|=\frac{n}{2}-1$. The graph $G$ is described by
\begin{itemize}
\item $G[A_i,B_i]$ is a complete bipartite graph for $i=1,2$, and
\item $G[A_1,B_2]\cong G[B_1,A_2]\cong Q(n/2+1,s-1)$.
\end{itemize}

We have $\delta(G)=\frac{n}{2}+s-2$. The fact that there exists no
$K_{s,t}$-factor is implied immediately by the fact that there is no
subgraph isomorphic to $K_{s,t}$ whose vertices would touch both $A_1$ and
$B_2$, or $A_2$ and $B_1$.

\subsection{Case $k$ is odd, $2s+1\geq t>s+1$}
Let $k=2l+1$, $n=k(s+t)$. Note that $\frac{n-t+s+2}{2}$ is an integer. Partition $A=A_1+A_2+A_*$, $B=B_1+B_2+B_*$,
$|A_1|=|A_2|=|B_1|=|B_2|=\frac{n-t+s+2}{2}$, $|A_*|=|B_*|=t-s-2$. The
graph $G$ is described by
\begin{itemize}
\item $G[A_i,B_i]$ is a complete bipartite graph for $i=1,2$,
\item $G[A_*,B_i]$ and $G[B_*,A_i]$ are complete bipartite graphs for $i=1,2$,
\item $G[A_1,B_2]\cong G[A_2,B_1]\cong P(\frac{n-t+s+2}{2},s-1)$,
\item the graph $G[A_*,B_*]$ is empty.
\end{itemize}
We have $\delta(G)=\frac{n+t+s}{2}-2$. To see that $G$ does not have
a $K_{s,t}$-factor, we argue as follows. Suppose for contradiction that $G$ has a
$K_{s,t}$-factor. Fix a $K_{s,t}$-factor of $G$. First, observe that there cannot be a copy
isomorphic to $K_{s,t}$ intersecting both $A_1\cup B_1$ and $A_2 \cup
B_2$. Let $r_1$ and $r_2$ be the number of copies of $K_{s,t}$
in the tiling whose color class of size $t$ touches $A_1$ and $B_1$,
respectively. Let $A_c$ and $B_c$ be vertices covered by these
$r_1+r_2$ copies. It holds
\begin{equation}\label{eq_subsets}
  A_1\subset A_c\subset A_1 \cup A_* \quad \mbox{and}\quad
 B_1\subset B_c\subset B_1 \cup B_*\,.
\end{equation}
 If $r_1\not = r_2$ then $\left||A_c|-|B_c|\right|\geq
t-s$, which contradicts~(\ref{eq_subsets}). Thus, $r_1 = r_2$. We conclude that
$$\frac{l(s+t)+s+1}{s+t}\leq r_1 \leq \frac{l(s+t)+t-1}{s+t} \mbox{,}$$
a contradiction to the integrality of $r_1$.
\subsection{Case $k$ is odd, $t=s+1$}
By $R(m,q)$ we denote (any of possibly many) bipartite graph
$R(m,q)=(R_1,R_2; E_R)$ with the following properties:
\begin{itemize}
\item $|R_1|=m, |R_2|=m-1$,
\item $R(m,q)$ does not contain any $K_{2,2}$,
\item for any vertex $x$ in $R_1$, it holds $\deg(x)\in\{q-1,q\}$, and
\item for any vertex $y$ in $R_2$, it holds $\deg(y)=q$.
\end{itemize}
For fixed $q$ and large $m$ the existence of such a graph $R(m,q)$ follows by a construction
analogous to the construction of the graph $Q(m,q)$.

Let $k=2l+1$. Partition $A=A_1+A_2$, $B=B_1+B_2$, $|A_1|=|B_1|=l(s+t)+s$,
$|A_2|=|B_2|=l(s+t)+s+1$. The graph $G$ is described by
\begin{itemize}
\item $G[A_i,B_i]$ is a complete bipartite graph for $i=1,2$,
\item $G[B_2,A_1]\cong G[A_2,B_1]\cong R((n+1)/2,s-1)$.
\end{itemize}
One immediately sees that $\delta(G)=\frac{n+t+s}{2}-2$ and no
$K_{s,t}$-tiling of $G$ exists.

\section{Upper bound}
We prove Theorem~\ref{thm_upper} in this section. 
The proof of Theorem~\ref{thm_upper} utilizes the previous work of Zhao~\cite{ZhaoBipTil}.
We will need the following lemma, which allows us
to find many vertex disjoint copies of certain stars. For
$h\in\mathbb{N}$, an $h$-star is a graph $K_{1,h}$, its
{\em center} is the unique vertex in the part of size one.
Moreover, for a graph $G$ and two disjoint sets $A,B\subset
V(G)$ we define
$$
\delta(A,B)=\min\{\deg(v,B)\::\:v\in A\}\;,\quad
\Delta(A,B)=\max\{\deg(v,B)\::\:v\in
A\}
$$
and 
$$
d(A,B)=\frac{e(A,B)}{|A||B|}\,.
$$

\begin{Lemma}[Zhao,~\cite{ZhaoBipTil}]\label{lemma12}
Let $1\le h\le \delta\le M$ and $0<c<1/(6h+7)$. Suppose
that $H=(U_1,U_2; E_H)$ is a bipartite graph such that
$||U_i|-M|\le cM$ for $i=1,2$. If
$\delta=\delta(U_1,U_2)\le cM$ and
$\Delta=\Delta(V_2,V_1)\le cM$, then we can find a family of vertex-disjoint $h$-stars, $2(\delta-h+1)$ of which have
centers in $U_1$ and $2(\delta-h+1)$ of which have centers
in $U_2$.
\end{Lemma}

As in~\cite{ZhaoBipTil} we distinguish between an extremal and a non-extremal case.
If we find a $K_{s+t,s+t}$-factor in $G$ we are done,
as each copy of $K_{s+t,s+t}$ can be split into two copies of
$K_{s,t}$ and hence we have a $K_{s,t}$-factor. Thus the theorem
stated next is just a corollary of~\cite[Theorem~4.1]{ZhaoBipTil}.
\begin{Theorem}[Zhao,~\cite{ZhaoBipTil}]
For every $\alpha>0$ and positive integers $s<t$, there exist $\beta>0$
and a positive integer $k_0$ such that the following holds for all
$n=k(s+t)$ with $k>k_0$. Given a bipartite graph $G=(A, B; E)$ with
$|A|=|B|=n$, if $\delta(G)> (\frac12-\beta)n$, then either $G$ contains a
$K_{s,t}$-factor, or there exist
\begin{equation*}
 A_1\subset A, \quad B_1\subset B \quad \mbox{such that}
 \quad |A_1|=|B_1|= \lfloor n/2 \rfloor, \quad d(A_1, B_1)< \alpha.
\end{equation*}
\end{Theorem}
Therefore, we reduce the problem to the extremal case. Let
$\alpha=\alpha(t)>0$ be small. As in the proof of Theorem~11 in~\cite{ZhaoBipTil},
define
\begin{equation*}
\begin{array}{r@{=}l r@{=}l}
  A'_1& \left\{x\in A: \deg(x, B_1)< \alpha^{\frac13}\,\frac{n}2\right\},
& B'_1& \left\{x\in B: \deg(x, A_1)< \alpha^{\frac13} \,
\frac{n}2\right\},
 \\
A'_2& \left\{x\in A: \deg(x, B_1)> (1-\alpha^{\frac13}) \frac{n}2\right\},
&B'_2& \left\{x\in B: \deg(x, A_1)> (1-\alpha^{\frac13})
\frac{n}2\right\},
\\
 A_0&A-A'_1-A'_2, & B_0&B-B'_1-B'_2,
 \\
 G_1&G[A'_1,B'_1], & G_2&G[A'_2,B'_2].
\end{array}
\end{equation*}
Similarly as in the proof of Theorem~11 in~\cite{ZhaoBipTil}, we assume that
removing any edge from $G$ would violate the minimum-degree
condition and then change $A'_i$ and $B'_i$ a little so that
$\Delta(G_1),\Delta(G_2)<\alpha^{\frac{1}{9}}n$. Vertices in $A_0\cup B_0$ are called {\em special}.

\subsection{$k$ is even}
To exhibit the existence of a tiling in this case, it is sufficient
to translate carefully the proof of Case~I of Theorem~11 from~\cite{ZhaoBipTil}.
We give a sketch of the proof below and refer the reader to
the corresponding places in~\cite{ZhaoBipTil} for more details.

Set $\mathcal{V}=(A'_1,B'_1,A'_2,B'_2)$. First assume, that no
member of $\mathcal{V}$ contains more than $n/2$ vertices. We add
vertices from $A_0$ and $B_0$ into sets of $\mathcal{V}$ in such a
way, that every set has size exactly $n/2$. Then, we may apply
arguments used in~\cite{ZhaoBipTil}, based on Hall's Marriage Theorem, to find a
$K_{s+t,s+t}$ tiling.

Next, assume that there is only one set in $\mathcal{V}$ which has more
than $n/2$ elements. Without loss of generality, assume that it is
$A'_2$, i.e., $|A'_2|=c>n/2$. Lemma~\ref{lemma12} applied to the graph
$G[A'_2,B'_2]$ yields the existence of $c-n/2$ disjoint
$s$-stars with centers in $A'_2$. We move the centers of the stars
into $A'_1$ and extend each of the stars into a copy of $K_{s,t}$ (each of these copies lies entirely in $A'_1\cup B'_2$, with the color class of size $s$ being contained in $B'_2$).
We distribute vertices of $B_0$ into $B'_1$ and $B'_2$ so, that
$|B'_1|=|B'_2|=n/2$. Then, it is easy to finish the entire tiling.
This is done in three steps. In the first step, we find in an
arbitrary manner $c-n/2$ copies of $K_{s,t}$ (disjoint with the previous
ones) in $G[A'_1,B'_2]$ placed in such a way, that the color-class
of size $s$ lies in $A'_1$. This
step ensures us, that the cardinalities of untiled (i.e., those vertices which are not covered by the partial $K_{s,t}$-factor) vertices in
the both color-classes of $G[A'_1,B'_2]$ are equal and divisible by
$s+t$. In the second step, all yet untiled vertices of
$G[A'_1,B'_2]$ which were originally special vertices are tiled. In
the third step, the tiling is in an analogous manner defined for
$G[A'_2,B'_1]$.

Now, assume that two diagonal sets of $\mathcal{V}$, say $A'_2$ and
$B'_1$ have sizes more than $n/2$. Then we apply separately Lemma~\ref{lemma12}
 to $G[A'_2,B'_2]$ and $G[A'_1,B'_1]$ to obtain families
$\mathcal{S}_A$ and $\mathcal{S}_B$ of disjoint $s$-stars with
centers in $A'_2$ and $B'_1$, such that
$|A'_2|-|\mathcal{S}_A|=|B'_1|-|\mathcal{S}_B|=n/2$. We move the
centers of the stars to $A'_1$ and $B'_2$ and proceed as in the
previous case.

The remaining case is when two non-diagonal sets from $\mathcal{V}$
have size more than $n/2$. Assume these are $A'_2$ and $B'_1$. We
apply Lemma~\ref{lemma12} to the graph $G[A'_2,B'_2]$ to obtain families
$\mathcal{S}_A,\mathcal{S}_B$ of disjoint $s$-stars with centers in
$A'_2$ and $B'_2$, such that
$|A'_2|-|\mathcal{S}_A|=|B'_2|-|\mathcal{S}_B|=n/2$. We proceed as
in the previous cases.

\subsection{$k$ is odd}
Let $k=2l+1$. We say that a set of special vertices ($A_0$ and/or $B_0$) is {\em small} if its size is less than $t-s$. Otherwise, it is called {\em big}.

We distinguish four cases.

\begin{itemize}
\item \emph{Both $A_0$ and $B_0$ are small.} Then there exist
$i,j\in\{1,2\}$, such that $|A'_i|, |B'_j|\ge l(s+t)+s+1$. If
$i=j$, then we apply Lemma~\ref{lemma12} to the graph $G_i$ and find
families $\mathcal{S}_A$, $\mathcal{S}_B$ of pairwise disjoint
$s$-stars with centers in $A'_i$ and $B'_i$ respectively, so that
$|A'_i|-|\mathcal{S}_A|=|B'_i|-|\mathcal{S}_B|=l(s+t)+s$. Move the
centers of the stars in $A'_{3-i}$ and $B'_{3-i}$. After the changes
we shall tile two graphs: $G[A'_1,B'_2]$ and $G[A'_2,B'_1]$. Note,
that both those graphs are not balanced. The tiling procedure is
analogous to the previous cases (when $k$ is even); the only difference is that one
copy of $K_{s,t}$ has to be found in the graphs first to make each
of them balanced.

If $i\not =j$, we can assume that $|A'_j|,|B'_i|\le
l(s+t)+s$. Since if this does not hold, then we could change one index and
continue as in the case when $i=j$. We will show
that one can add vertices to $A'_j$ and to $B'_i$ so that
$|A'_j|=l(s+t)+s$ and $|B'_i|=l(s+t)+t$. Then, the existence of the
tiling will follow by standard arguments. We apply Lemma~\ref{lemma12} to the
graph $G_j$ to obtain a family of $|B'_j|-(l(s+t)+s)$ vertex
disjoint $s$-stars with centers in $B'_j$ and end-vertices in
$A'_j$. If we moved all the centers to $B'_i$ and all the vertices
of $B_0$, the cardinality of $B'_i$ would be
$$|B'_i|+(|B'_j|-(l(s+t)+s))+|B_0|=l(s+t)+t \; .$$
The same applies for $A'_j$. Therefore, by removing some of the
vertices, we may attain $|A'_j|=l(s+t)+s$ and $|B'_i|=l(s+t)+t$.
Then, the existence of a tiling follows.

\item \emph{$A_0$ is small and $B_0$ is big.} Then at least one $B'_i$
(say $B'_2$) has at most $l(s+t)+s$ vertices. Lemma~\ref{lemma12} asserts that
we can find a family $\mathcal{S}_B$ of disjoint $s$-stars with
centers in $B'_1$ and end-vertices in $A'_1$, such that
$|B'_1|-|\mathcal{S}_B|\le l(s+t)+s$. This implies, that we can find
vertices (in $B_0$ or centers of the stars of $\mathcal{S}_B$) which
can be moved to $B'_2$ so that $|B'_2|=l(s+t)+t$.

As $A_0$ is small, one of $A'_1$ and $A'_2$ must have at least
$l(s+t)+s+1$ vertices. The tiling can be found by standard arguments
if we achieve to have $|A'_1|=l(s+t)+s$. If $|A'_1|>l(s+t)+s$, Lemma~\ref{lemma12} yields existence of a family $\mathcal{S}_A$ of disjoint
$s$-stars with centers in $A'_1$ and end-vertices in $B'_1$ such
that $|A'_1|-|\mathcal{S}_A|=l(s+t)+s$. Moving the centers to $A'_2$,
we achieve $|A'_1|=l(s+t)+s$. Assume that $|A'_1|\le l(s+t)+s$. The
size of $A'_2$ is $k(s+t)-|A'_1|-|A_0|>l(s+t)+s$. Lemma~\ref{lemma12} yields
existence of a family $\mathcal{S}_A$ of disjoint $s$-stars in $G_2$
centered in $A'_2$ with the property that
$|A'_1|+|\mathcal{S}_A|=l(s+t)+s$. Moving the centers to $A'_1$
yields demanded $A'_1=l(s+t)+s$.

\item \emph{$A_0$ is big and $B_0$ is small.} The analysis of this case is analogous to
the previous one.

\item \emph{Both $A_0$ and $B_0$ are big.} We shall show in the next paragraph, that we
can achieve $A'_1$ to be of size $l(s+t)+s$ and of size $l(s+t)+t$.
An analogous procedure can be used to show the same property for the
set $B'_1$. Then, the existence of the tiling follows immediately;
one prescribes the cardinalities of $A'_1$ and $B'_1$ to be equal to
the same number $l(s+t)+s$.

If $|A'_i\cup A_0|<l(s+t)+t$ for some $i\in\{1,2\}$, then we have
$|A'_{3-i}|>l(s+t)+s$. Appealing to Lemma~\ref{lemma12}
we can remove centers of $s$-stars from $A'_{3-i}$  in
such a way that $|A'_{3-i}|=l(s+t)+s$. Also, by moving $t-s$ vertices from the big set $A_0$ to $A'_{3-i}$ arrive at $|A'_{3-i}|=l(s+t)+t$. Then, the partial $K_{s,t}$-tiling can be extended to a $K_{s,t}$-factor.

Finally, if both $|A'_1|\le l(s+t)+s$ and $|A'_2|\le l(s+t)+s$ then we redistribute some vertices (again, appealing to Lemma~\ref{lemma12}, and using the set $A_0$) to arrive at the situation when $|A'_1|=l(s+t)+s$, $|A'_2|=l(s+t)+t$. Then the tiling can be extended as before.
\end{itemize}

\subsection*{Acknowledgement}
We thank a careful referee for suggesting several
improvements in the presentation.
\bibliographystyle{plain}
\bibliography{bibl}
\end{document}